\documentclass[11pt]{amsart}

\usepackage{tikz}
\usepackage{amsmath}
\usepackage{amsthm}
\usetikzlibrary{calc}
\usepackage{amsfonts}
\usepackage{amssymb}
\usepackage[all]{xy}
\usepackage{alltt}
\usepackage{enumitem}
\usepackage{stmaryrd}
\usepackage{comment}
\usepackage{xspace}
\usepackage{comment}
\usepackage{libertine}
\usepackage[libertine]{newtxmath}
\usepackage{graphicx}
\theoremstyle{plain}
\newtheorem{proposition}{Proposition}[section]

\theoremstyle{definition}

\newtheorem{remark}[proposition]{Remark}

\usepackage[margin=1.25in]{geometry}

\usepackage[
        colorlinks=true,urlcolor=blue,linkcolor=blue,citecolor=blue
]{hyperref}

\begin{document}

% Article information
\title{The complexities of falling freely}
\date{\today}
 \author{Anindya Sen} 
\address {University of Otago, NZ}
\email{anindya.sen@otago.ac.nz}

% Author information
\author{Sunil Chebolu}
\address{Department of Mathematics \\
Illinois State University \\
Normal, IL 61790, USA}
\email{schebol@ilstu.edu}

\begin{abstract}

Suppose you drop a coin from 10 feet above the ground. How long does it take to reach the ground? This routine exercise is well-known to every AP physics and calculus student: the answer is given by a formula that assumes constant acceleration due to gravity. But what if you ask the same question in the more realistic scenario of non-constant acceleration following an inverse square law? In this article, we explain the analysis of this realistic scenario using freshman-level calculus and examine some implications. As a bonus, we also answer the following intriguing question: Suppose the Earth were to instantaneously collapse to a mathematical point at its center. How long would it take for us surface dwellers to fall to the center? 
\end{abstract}

\maketitle
\thispagestyle{empty}

%Updated and Revised CMJ Template 12.20.2023 by Bonnie Ponce and Bev Ruedi copyright the Mathematical Association of America and all rights reserved.

% PLEASE % COMMENT OUT AUTHOR BIOS AND DO NOT ADD AUTHOR PHOTO(S) UNTIL YOU RECEIVE YOR PROVISIONAL ACCEPT LETTER.

\begin{comment}
\begin{biog} %comment out in initial submission to allow for double blind reviewing
\item[\biogpic{\includegraphics[width=84pt]{anindya.jpeg}}Anindya Sen] (anindya.sen@otago.ac.nz) received his Ph.D. in Mathematics from the University of Chicago in 2005. After a stint as an investment banker at Goldman Sachs, he is currently an academic at the University of Otago in New Zealand. Anindya recommends appreciating the joy of falling freely by skydiving.

\vspace{2 cm} %this space may not be needed if the author's biography is long enough.

\item[\biogpic{\includegraphics[width=84pt]{sunil.jpg}}Sunil Chebolu] (schebol@ilstu.edu) received his Ph.D. in Mathematics from the University of Washington in 2005. He is currently a professor at Illinois State University, where he takes great pleasure in working with students. In his free time, Sunil enjoys the thrill of trail running and the wonder of exploring the cosmos with his telescopes.
\end{biog}

%\vspace{2 cm} %this space may not be needed if the author's biography is long enough
\end{comment}

\section{Introduction}
Let us begin with a standard problem from high school physics. Consider an object falling freely from rest under gravity (typically assumed to be Earth's
gravity and ignoring air resistance). How long does it take to fall by a height $h$? We will let $T$ denote this time.

A standard solution to find a formula for $T$ goes as follows. 
Let $g$\,($\,\approx 9.8 \, m/s^2$ for Earth) be the acceleration due to gravity, and let $s(t)$ denote the distance fallen in the first $t$ seconds.   The acceleration is constant since the object is in free fall, meaning  $d^2 s/dt^2= g$. Integrating this equation twice and applying the initial conditions ($s(0)=0$ and $s'(0)=0)$ gives $s(t) = \frac{1}{2}gt^2  $. Setting this equation equal to $h$ and solving for $T$ then gives the high-school formula:
\begin{equation} \label{simpleformula}
\boxed{
T = \sqrt{\frac{2h}{g}}.}
\end{equation}

This formula, of course, assumes that gravitational acceleration is constant.
In reality, gravitational acceleration follows an inverse square law. How does the above formula change when we consider this more realistic scenario? To explain this scenario, we must recall two important Newtonian laws of physics.

A fundamental law of motion states that when a net force of magnitude $F$ acts on a mass $m$, it produces an acceleration with magnitude $a$ such that 
\[ F = ma.\] 

The second is one of the most famous laws in physics and astronomy: Newton's universal law of gravity. It states that the magnitude $F$ of the gravitational force of attraction between objects of masses $M$ and $m$ that are separated by a distance $r$ is given by 
\[F = \frac{G M m}{r^2}, \]
where $G$ is the universal constant of gravity whose value is approximately equal to 
	$6.67430 \times 10^{-11} \ \  m^3 kg^{-1} s^{-2}$.
%These two equations are incredibly powerful; to some extent, they explain the motion of all objects in the universe, from falling apples to orbiting planets and satellites!

Equating these two values of $F$ and canceling out $m$ on both sides gives the  inverse square law:  the gravitation acceleration of an object that is at a distance $r$ units from a point mass $M$ (mass of Earth in the case of free fall) is given by
\begin{equation} \label{g-G}
g(r) = \frac{GM}{r^2}.
\end{equation}
(Note that the acceleration due to gravity does not depend on the mass $(m)$ of the body being attracted. This is a fact with deep-reaching implications beyond the scope of this article.)

So now we can ask a more advanced question.
Suppose we have an object falling freely from rest under the gravity of a point mass $M$, where the resulting acceleration is given by equation (\ref{g-G}). Suppose its initial distance from the point mass is $R_0$.
How long does it take to fall to a distance $R_1$, where $0 \le R_1 \le R_0$?

In this case, the answer is given by the rather more formidable formula:

\begin{equation} \label{grandequation}
\boxed{
T = \frac{R_0^{3/2}}{\sqrt{2GM}} \left( \frac{\pi}{2} - \sin^{-1}\left(\sqrt{\frac{R_1}{R_0}} \right)  + \sqrt{\frac{R_1}{R_0} \left(1-\frac{R_1}{R_0}\right) }\right).
}
\end{equation}

Several questions quickly spring to mind.
\begin{enumerate}
    \item[Q1.] How did we arrive at this intimidating formula (\ref{grandequation})?
    \item[Q2.] 
How could it possibly connect to the much simpler formula (\ref{simpleformula})?
\item[Q3.] It is also surprising that the number $\pi$ shows up in the formula (\ref{grandequation}). This is usually a sign that
circles, spheres, or something similar are involved. How could that be?

\end{enumerate}

We shall answer these questions in this article using only undergraduate-level calculus. All calculus results used here can be found in any standard calculus textbook; see \cite{stewart} for instance.

\section{Formula Derivation} In this section, we will derive formula (\ref{grandequation}). Let us begin with a recap of the framework and notation. Consider the gravitational field induced by a fixed point mass $M$. An object of mass $m$ at rest in this field is released at a distance $R_0$ from the point mass at time $t=0$. Let $r = r(t)$ denote the distance between the falling object and the point mass at time $t$, and let $T$ denote the time the object takes to reach a point whose distance from the point mass is $R_1$. 

Using the inverse square law of gravitational acceleration from the previous section, we get the following initial value problem:
\[ \frac{d^2r}{dt^2} = - \frac{GM}{r^2}, \; \;  r(0) = R_0, \; \; r'(0) = 0.\]
Even though this is a second-order non-linear differential equation connecting $r$ and $t$, it can be solved using elementary methods as follows.

We multiply both sides of the above differential equation with $dr/dt$ and then integrate. This gives (using Newtonian notation of derivatives): 
\begin{eqnarray*}
r''& = &  - \frac{GM}{r^2} \\
r' r'' &= &  - r' \frac{GM}{r^2} \\
\int r' r'' \, dt & = &  \int - \frac{GM}{r^2} r'\, dt, \\
\frac{r'^2}{2} & = &   \frac{GM}{r} + C \;\;  \text{for some constant $C$}\\
r'^2 & = &   \frac{2GM}{r} + C'  \;\;  \text{for some constant $C'$.}
\end{eqnarray*}

Setting $t=0$ in the last equation and invoking the given initial conditions, we see that $C' = -2GM/R_0$. Also, note that $r$ can only decrease under a free fall, so $r' < 0$. This gives 
\begin{equation}
    r'(t) = \frac{dr}{dt}= - \sqrt{2GM\left( \frac{1}{r} - \frac{1}{R_0}\right)}, \label{r'(t)}
\end{equation}
with the initial condition $r(0) = R_0$. \\
For future reference, we write $R_1 = KR_0$ for some $K \in [0, 1]$.  
We now want to find $T$ such that $r(T) = R_1$, for some $0 \le R_1 < R_0$. 
To this end, we begin with our  (separable) differential equation and integrate both sides after separating the variables in question. This gives the following equations:

\begin{eqnarray*}
 %\frac{dr}{dt} & = &  - \sqrt{2GM\left( \frac{1}{r} - \frac{1}{R_0}\right)  \\
    dt & = & - \frac{dr}{ \sqrt{2GM\left( \frac{1}{r} - \frac{1}{R_0}\right) }} \\
    \int_0^T \, dt & = & \int_{R_0}^{R_1} - \frac{dr}{ \sqrt{2GM\left( \frac{1}{r} - \frac{1}{R_0}\right) }}\\
    T  & = & \int_{R_1}^{R_0}  \frac{dr}{ \sqrt{2GM\left( \frac{1}{r} - \frac{1}{R_0}\right) }}  = \frac{1}{\sqrt{2GM}} \int_{R_1}^{R_0} \sqrt{\frac{R_0r}{R_0 -r}} \, dr \\
    T & = & \sqrt{\frac{R_0}{2GM}} \int_{R_1}^{R_0} \frac{\sqrt{r}}{\sqrt{R_0-r}} \; dr
\end{eqnarray*}
Evaluating the integral on the right-hand side of the last equation is a standard exercise using a trigonometric substitution: $r = R_0 \sin^2 \theta$. (Note that this is a valid substitution because $r$ is positive and bounded above by $R_0$.)  Note that, under this change of variable $r \rightarrow \theta$, we have  $R_0 \mapsto \pi/2$ and $R_1 \mapsto \sin^{-1} \sqrt{\frac{R_1}{R_0}} = \sin^{-1}\sqrt{K}$. 

Substituting these in the last integral, we get:
\begin{eqnarray*}
  %  T & = &   \sqrt{\frac{R_0}{2GM}} \int_{R_1}^{R_0} \frac{\sqrt{r}}{\sqrt{R_0-r}} \; dr \\
   T  & = &  \sqrt{\frac{R_0}{2GM}} \int_{\sin^{-1}\sqrt{K}}^{\pi/2} (\sqrt{R_0}\sin \theta \, 2 R_0 \sin \theta \cos \theta)/\sqrt{R_0} \cos \theta  \, d\theta\\
    & = & \sqrt{\frac{R_0}{2GM}} 2R_0 \int_{\sin^{-1}\sqrt{K}}^{\pi/2} \, \sin^2 \theta \, d\theta \\
    & = & \sqrt{\frac{R_0}{2GM}} 2R_0  \frac{1}{2} (\theta - \sin \theta \cos \theta) \bigg \rvert_{\sin^{-1}\sqrt{K}}^{\pi/2}. 
\end{eqnarray*}   

%Now, we know that 
%\[\sin^2 \theta = \frac{1}{2}(1 -\cos 2\theta) \implies \int \, \sin^2 \, d\theta = \frac{1}{2} (\theta - \frac{1}{2} \sin 2 \theta) =  \frac{1}{2} (\theta - \sin \theta \cos \theta).\]

%Plugging this into our integral, we get 
%\[T =  \sqrt{\frac{R_0}{2GM}} 2R_0  \frac{1}{2} (\theta - \sin \theta \cos \theta) \bigg \rvert_{\sin^{-1}\sqrt{K}}^{\pi/2}. \]
%It is a standard exercises using the double-angle formula ($\cos 2x = 1 - 2\sin^2 x$)  to show that $\int \, \sin^2 \, d\theta = \frac{1}{2} (\theta - \sin \theta \cos \theta) $
Finally, note that 
$\theta = \sin^{-1} \sqrt{K} \implies \cos \theta = \sqrt{1 - \sin^2 \theta} = \sqrt{1-K}.$
Putting all these together and doing some algebra, we  arrive at:
\[
T = \frac{R_0^{3/2}}{\sqrt{2GM}} \left( \frac{\pi}{2} - \sin^{-1}(\sqrt{K})  + \sqrt{K (1- K) }\right).
\]
Finally, substituting $K = R_1/R_0$, we get the exact formula mentioned in the introduction.

\vskip 5mm
\begin{remark} 
In the above formula, it appears as though the falling time \( T \)  is proportional to \( R_0^{3/2} \) for a fixed point mass $M$. This might seem problematic since \( R_0 \) depends on the unit of length. For example, if \( R_0 \) is measured in centimeters instead of meters, \( R_0 \) increases by a factor of 100, and \( R_0^{3/2} \) increases by a factor of 1000. This might suggest that \( T \) would also increase by a factor of 1000, but clearly, this change in units should have no effect on the value of the falling time, $T$.

The resolution comes from the dimensional analysis of factor \( \frac{R_0^{3/2}}{\sqrt{2GM}} \). (The term inside the brackets is a function of $K (= R_1/R_0)$ and therefore dimensionless and unaffected by unit changes.)
 Note that \( R_0 \) has units of length,  \( M \) has units of mass, and finally \( G (= \frac{gR_0^2}{M}, \) see Equation (\ref{g-G})) has units of \( \text{length}^3 \, \text{mass}^{-1} \, \text{time}^{-2} \).
Therefore, if the unit of length changes by a factor of \( \lambda \), then \( R_0 \) scales as \( R_0 \rightarrow \lambda R_0 \) and \( G \) scales as \( G \rightarrow \lambda^3 G \). As a result, \( \frac{R_0^{3/2}}{\sqrt{2GM}} \) remains unchanged.

A similar analysis shows that the factor \( \frac{R_0^{3/2}}{\sqrt{2GM}} \)  is also invariant under a change of the unit of mass.
In fact, its value only depends on the unit of time measurement, which is precisely what we should expect.

\end{remark}

\section{Reconciliation}
As promised, we derived the exact formula for free fall time under gravity that follows the inverse square law. However, it would be reassuring to reconcile it with the simpler formula from school days! Among other things, this would reinforce our confidence that we haven't inadvertently made a mistake in the derivation.
%– which is quite likely when the math gets complicated. 
Let us begin by recording the two formulas one below the other for a quick comparison:
\begin{eqnarray*}
  \text{High School Formula:} \; \; \;  T & = &   \sqrt{\frac{2h}{g}}. \\
\text{Exact Formula:} \; \; \;  T & = & \frac{R_0^{3/2}}{\sqrt{2GM}} \left( \frac{\pi}{2} - \sin^{-1}(\sqrt{K})  + \sqrt{K (1- K) }\right).
\end{eqnarray*}   

Where do we begin? Let us first compare the variables that appear in both formulas. Since $h = R_0 - R_1$ and $K=R_1/R_0$, both formulas involve $R_0$ and $R_1$. 
However, the high school formula involves $g$, whereas the exact formula involves $G$. So, we need an expression for the constant gravitational acceleration $g$ in terms of $G$. Since we are assuming that the entire free fall is taking place at a distance very close to $R_0$ from the point mass $M$, using Equation (\ref{g-G}), we can write  
\[g = \frac{GM}{R_0^2}.\]

Let us plug this into the high school formula and simplify it until it starts to resemble the term outside the brackets in the exact formula. This gives:
\[
 T  =    \sqrt{\frac{2h}{g}} 
      =  \frac{\sqrt{2h}}{\sqrt{\frac{GM}{R_0^2}}}
      = \frac{R_0 \sqrt{2h}}{\sqrt{GM}}
= \frac{R_0 2 \sqrt{h}}{\sqrt{2GM}}= \frac{R_0^{3/2}}{\sqrt{2GM}} 2 \sqrt{\frac{h}{R_0}}
.\]

Once we  get rid of the common factor $R_0^{3/2}/\sqrt{2GM}$, the high school and exact formulas boil down to comparing the values:
\[2 \sqrt{\frac{h}{R_0}} \; \; \; \text{ versus } \; \; \;\frac{\pi}{2} - \sin^{-1}(\sqrt{K})  + \sqrt{K (1- K) }.\]

Also note that $K = R_1/R_0 = (R_0-h)/R_0 = 1 - h/R_0$. This implies that $1 - K = h/R_0$. Setting $\epsilon = h/R_0$, the two expressions under consideration are 
\[ 2\sqrt{\epsilon} \; \; \; \text{ versus } \; \; \; \frac{\pi}{2} - \sin^{-1}(\sqrt{1 -\epsilon})  + \sqrt{\epsilon (1- \epsilon) }.\] \label{twothings}

In the limit, where $\epsilon = h/R_0$ is extremely small, could these two expressions be the same? (Note that since the functions involved here are not differentiable at the origin, standard tangent line approximations are not feasible.) 
%It is straightforward to see that:
%\[ \lim_{\epsilon \rightarrow 0} \;  \frac{\pi}{2} - \sin^{-1}(\sqrt{1 -\epsilon})  + \sqrt{\epsilon (1- \epsilon) } - 2\sqrt{\epsilon} = 0. \]
 When $\epsilon \approx 0$, $\sqrt{1- \epsilon} \approx 1$. This means that  $\sqrt{\epsilon (1- \epsilon)} \approx \sqrt{\epsilon}$.
 This approximation takes into account $\sqrt{\epsilon}$. Now we need to compare:
 
 %Subtracting $\sqrt{\epsilon}$ from both the expressions, we are left with 
\[ \sqrt{\epsilon} \; \; \; \text{ versus } \; \; \; \frac{\pi}{2} - \sin^{-1}(\sqrt{1 -\epsilon}).   \]
We will be done if we can show that these two quantities are approximately equal when $\epsilon$ is close to zero.
 To see this, consider a right triangle with a hypotenuse of 1,  where the side opposite to angle $\theta$ has length $\sqrt{\epsilon}$. Consequently, the adjacent side will have length $\sqrt{1 -\epsilon}$.

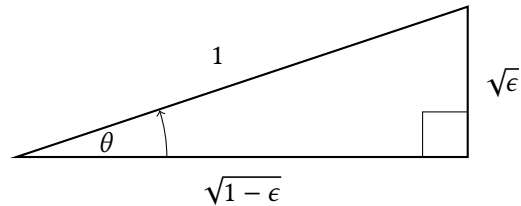
\begin{center} 
\begin{tikzpicture}[scale=2]
    % Define the coordinates for the vertices of the triangle
    \coordinate (A) at (0,0);
    \coordinate (B) at (3,0); % Adjust this to change the length of side 'b'
    \coordinate (C) at (3,1); % Adjust this to make the angle skinnier

    % Draw the triangle
    \draw[thick] (A) -- (B) -- (C) -- cycle;

    % Label the sides just outside the triangle, extremely close to the midpoints
    \node[below] at ($(A)!0.5!(B) + (0,-0.05)$) {$\sqrt{1 - \epsilon}$};
    \node[right] at ($(B)!0.5!(C) + (0.05,0)$) {$\sqrt{\epsilon}$};
    \node[above left] at ($(A)!0.5!(C) + (-0.05,0.05)$) {$1$};

    % Label the angle theta inside the triangle with a skinnier arc that meets the opposite side
    \draw[->] (1,0) arc[start angle=0,end angle=18,radius=1cm]; % Reduced the angle to make it skinnier
    \node at (0.6,0.1) {$\theta$};

    % Indicate the right angle
    \draw ($(B) + (0,0.3)$) -- ($(B) + (-0.3,0.3)$) -- ($(B) + (-0.3,0)$);
\end{tikzpicture}
\end{center}

Then note that 
\[ \frac{\pi}{2} - \sin^{-1}(\sqrt{1 -\epsilon})  =  \cos^{-1}(\sqrt{1 -\epsilon}) = \theta \approx \sin \theta =    \sqrt{\epsilon},\]
where the  approximation $(\sin \theta \approx \theta)$ is  the well-known first degree Taylor approximation for $\sin \theta$; note that $\theta \approx 0$ when $\epsilon \approx 0$. This shows the desired approximation.

 %T & =&  \frac{R_0^{3/2}}{\sqrt{2GM}} \left( \frac{\pi}{2} - \sin^{-1}\left(\sqrt{\frac{R_1}{R_0}} \right)  + \sqrt{\frac{R_1}{R_0} \left(1-\frac{R_1}{R_0}\right) }\right).
 
To summarize: When $h/R_0$ is very small,  both terms $\sqrt{K (1- K) }$ and $\frac{\pi}{2} - \sin^{-1}(\sqrt{K})$ of our exact formula are approximately equal to $\sqrt{h/R_0}$.  As a result, the expression $\frac{\pi}{2} - \sin^{-1}(\sqrt{K})  + \sqrt{K (1- K) }$ approximately equal to $2\sqrt{h/R_0}$. This approximation aligns our exact formula with the high school formula, thereby reconciling the two.

The high school formula works when the free-fall distance is very small compared to $R_0$. Let us consider the other extreme.
\emph{Suppose the Earth were to instantaneously collapse to a mathematical point at its 
center. How long would it take us surface dwellers to fall to the center?} Alternatively, if all the electrostatic forces of repulsion were magically turned off, how long
would it take the Earth to collapse to a point mass? (The latter situation is relevant when
analyzing the gravitational collapse of massive stars.)

Let us denote this time $T_C$ – for “Time to collapse.” With our exact free fall formula in hand, calculating $T_C$ is quite simple. We just set $R_0$ as the radius $R$ of the Earth and $R_1 = 0$. The latter implies $K =0$, and consequently, both $\sin^{-1} \sqrt{K}$  and $\sqrt{K(K-1)}$ become $0$, so we get:
\begin{equation} \label{collapse}
    T_C = \frac{\pi}{2} \frac{R^{3/2}}{\sqrt{2GM}}.
\end{equation}
Plugging in the radius $R$, the mass $M$ of the Earth, and the $G$ value in the above formula, we get about 15 minutes (or 896 seconds, to be more precise.)

\section{The $\pi$ factor}
As mentioned in the introduction, it is somewhat surprising that the constant $\pi$ enters our exact formula. In Section 2, we can see exactly what happens, but even then, it just “drops out of the
mathematics,” so to speak. Can we get a more intuitive understanding of what is happening
here?

A first clue comes from looking at realistic situations involving motion under gravity. Significant gravitational acceleration typically happens due to objects of planetary size or
larger, and the resulting motion takes place over thousands of kilometers or much more.
In such situations, vertical free fall, like the one we are considering, is rare – objects typically have
significant angular momentum with respect to the gravitating mass, and the result is some
kind of orbital motion.

The simplest case is when an object of mass $m$ orbits around a point mass $M$ in a perfect circle with radius $R$ and constant angular velocity $\omega$. Once again, this motion is governed by Newton's universal law of gravity: 
\begin{equation} \label{Newton}
    \frac{GMm}{R^2} = ma,
\end{equation} 
where $a$ is the acceleration of uniform circular motion. It can be shown using calculus that the acceleration of uniform circular motion is given by $a = R\omega^2$; see \cite{stewart}.
Plugging this value of $a$ in Equation \ref{Newton} and simplifying gives: 
\[ \omega^2R = \frac{GM}{R^2}  \implies \omega = \sqrt{\frac{GM}{R^3}}.\]
From there, we get a formula for the orbital period:
\[ T = \frac{2\pi}{\omega} = 2\pi \frac{R^{3/2}}{\sqrt{GM}}.\]

Note the similarity to the formula for $T_C$ in (\ref{collapse}). In fact, it is identical up to a constant factor of $4\sqrt{2}$.

However, surely this is just a coincidence?
What could the circular orbital period have to do with vertical free fall?  The situation gets more intriguing once we consider the most realistic case. Unlike vertical free fall and circular orbits, the most realistic situation in planetary dynamics
is elliptical orbits. In fact, according to Kepler's first law, planetary orbits are ellipses with the Sun at one focus. 

So, what is the time period for an elliptical orbit?
The orbital period $T_p$ is given by:
\begin{equation} \label{ellipticperiod}
 T_p = 2\pi \frac{R_*^{3/2}}{\sqrt{GM}}.
\end{equation}
Here, $R_* = \frac{1}{2} (R_{\max} +R_{\min})$, where 
$R_{\max}$ and $R_{\min}$ are the maximum and minimum distances between the two objects in question. (For instance, in the case of a lunar orbit around the Earth, these distances correspond to points known as apogee and perigee.) This formula is derived in any standard textbook on classical mechanics using the law of conservation of angular momentum; see \cite{anil}, for instance.

This formula is very similar to that of the orbital period for circular motion. But how does this help?  Well, consider an ellipse of extremely high eccentricity – very long and narrow. In other words, the minor axis has a length very close to 0, so it “almost” looks like a straight
line.

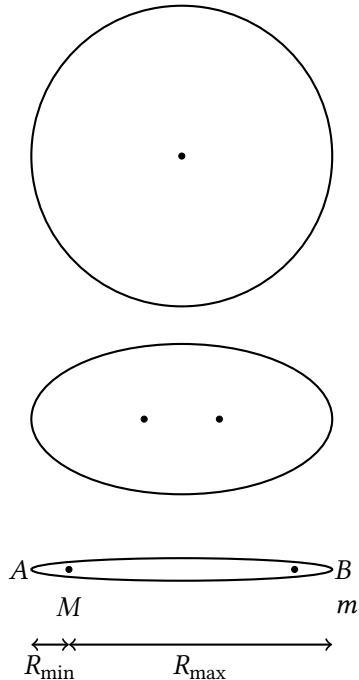
\begin{figure}[!h] 
\begin{center}
\begin{tikzpicture}
    % Circle
    \draw[thick] (0, 2) circle [radius=2];
    \filldraw (0, 2) circle [radius=0.04];
    
    % Ellipse
    \draw[thick] (0, -1.5) ellipse [x radius=2, y radius=1];
    \filldraw (-0.5, -1.5) circle [radius=0.04];
    \filldraw (0.5, -1.5) circle [radius=0.04];
    
    % Very skinny horizontal ellipse with points and labels
    \draw[thick] (0, -3.5) ellipse [x radius=2, y radius=0.15];
    \filldraw (-1.9, -3.5) circle [radius=0.00] node[anchor=east] {$A$};
    \filldraw (1.9, -3.5) circle [radius=0.00] node[anchor=west] {$B$};
    \filldraw (-1, -3.5) circle [radius=0.0] node[anchor=north] {};
    \filldraw (1.5, -3.5) circle [radius=0.04];
    \filldraw (-1.5, -3.5) circle [radius=0.04];
      % Label for M
    \node at (-1.5, -4) {$M$};
    
    % Label for m
    \node at (2.2, -4) {$m$};

    % Arrows and labels for distances
    \draw[thick][<->] (-2, -4.5) -- (-1.5, -4.5) node[midway, below] {$R_{\min}$};
    \draw[thick][<->] (-1.5, -4.5) -- (2, -4.5) node[midway, below] {$R_{\max}$};
\end{tikzpicture}

\end{center}
\caption{Ellipses with increasing eccentricity $e$. Circle on the top with $e=0$, and a skinny ellipse at the bottom with $e$ close to 1.}
\label{ellipses}

\end{figure}

Let $A$ and $B$ be the vertices of the ellipse - two ends of the major axis of the ellipse.
The point mass, $M$, will be at one focus – that can be assumed to be very close to $A$, and
suppose the falling object with mass $m$ is at $B$; see Figure \ref{ellipses}.

Do you see where this is going? The idea is that, in this case, \emph{the free-fall time, $T_C$ from $B$ to $M$, will be virtually identical to the time to follow the orbit from $B$ to $A$!}

Let us calculate this, then!
$R_{\max} =$ distance $MB = R_0$ and $R_{\min} \approx 0$ since the ellipse is really narrow. This gives:
\[ R_* = \frac{1}{2} (R_{\max} +R_{\min}) = R_* = \frac{R_0}{2}.\]

Plugging this into Equation \ref{ellipticperiod} gives us: 
\[ T_p = 2\pi \frac{R_*^{3/2}}{\sqrt{GM}} = 2\pi \frac{R_0^{3/2}}{\sqrt{GM}}\left(\frac{1}{2}\right)^{\frac{3}{2}} = \pi \frac{R_0^{\frac{3}{2}}}{\sqrt{2GM}}.\]
But wait, we are not interested in the full orbital period – we want the time taken to follow
the orbit from $B$ to $A$. This is \emph{half of} $T_P$!
Thus, finally, we get:
\[ T_C = \lim_{\text{eccentricity} \rightarrow 1} \frac{1}{2}T_p = \frac{\pi}{2} \frac{R^{3/2}}{\sqrt{2GM}}.\]
Our formula for free-fall time, which reconciled beautifully with the high school formula for
small distances, also matched up very neatly with formulas derived from the equations of
orbital dynamics for large distances.
The unexpected appearance of $\pi$ in the formula is now seen to be a consequence of the
connection between vertical free fall and the more general dynamics of orbital motion. Indeed, as anticipated, there are circles hidden behind the formula.
On this high note, we conclude the article.

%\begin{acknowledgment}
%The authors wish to thank the anonymous referees for their comments and suggestions.  
%\end{acknowledgment}

\begingroup
\raggedright

\bibliographystyle{plain}

\endgroup

%\begingroup
%\raggedright

%\bibliographystyle{alpha}
%\bibliography{citations}
%\endgroup

\end{document}